\documentclass[11pt]{article}
\usepackage[all]{xy}

\usepackage{epic}
\usepackage{amsmath}[1996/11/01]
\usepackage{amssymb,amsthm,amsfonts,latexsym,epsfig,graphics}
 \def\beql#1#2\eeql{\begin{equation}\label{#1}#2\end{equation}}

\advance \textheight by -1 \topmargin
\advance \textheight by -1 \headheight
\advance \textheight by -1 \headsep
\advance \textheight by -1 \footskip
\vsize = \textheight
\hsize = \textwidth
\newcommand{\knubbel}{\begin{picture}(2,2)(5,-3)  \put(-2,1){\circle*{4}} \end{picture}}
\newcommand{\operp}{\mathop{\bigcirc\!\!\!\!\!\!\!\perp\,}} 
\DeclareMathOperator{\Gal}{Gal}

\DeclareMathOperator{\GL}{GL}
\DeclareMathOperator{\Aut}{Aut}

\DeclareMathOperator{\im}{im}

\DeclareMathOperator{\Tr}{trace}
\DeclareMathOperator{\Cl}{Cl}

\newtheorem{theorem}{Theorem}[section]

\newtheorem{corollary}[theorem]{Corollary}

\newcommand{\bew}{\noindent\underline{Proof.}\ }
\newtheorem{rem}[theorem]{Remark}
\newtheorem{remark}[theorem]{Remark}
\newtheorem{lemma}[theorem]{Lemma}
\newtheorem{proposition}[theorem]{Proposition}
\newtheorem{definition}[theorem]{Definition}
\newtheorem{defn}[theorem]{Definition}

\newcommand{\Z}{{\mathbb{Z}}}
\newcommand{\Q}{{\mathbb{Q}}}
\newcommand{\F}{{\mathbb{F}}}

\newcommand{\R}{{\mathbb{R}}}
\newcommand{\E}{{\mathbb{E}}}
\newcommand{\A}{{\mathbb{A}}}
\newcommand{\D}{{\mathbb{D}}}

\newcommand{\trace}{\mbox{trace}}

\newcommand{\eb}{\phantom{zzz}\hfill{$\square $}\smallskip}

\begin{document}
\begin{center}
{\Large {\bf Automorphisms of extremal unimodular lattices in dimension 72.}}\\ 
\vspace{1.5\baselineskip}
{\em Gabriele Nebe} \\
\vspace*{1\baselineskip}
Lehrstuhl D f\"ur Mathematik, RWTH Aachen University\\
52056 Aachen, Germany \\
 nebe@math.rwth-aachen.de \\
\vspace{1.5\baselineskip}
\end{center}

{\sc Abstract.}
{\small 
The paper narrows down the possible automorphisms of extremal
even unimodular lattices of dimension 72. 
With extensive computations in {\sc Magma} using the very sophisticated algorithm 
for computing class groups of algebraic number fields written
by Steve Donnelly it is shown that the extremal even unimodular lattice
$\Gamma _{72}$ from \cite{dim72} is the unique 
extremal even unimodular lattice of dimension 72 that 
admits a large automorphism, where a $d\times d$ matrix 
is called large, if its minimal polynomial has an
irreducible factor of degree $> d/2$.
}

\section{Introduction}

This paper is part of a series of papers classifying extremal lattices 
with given automorphism.

A {\em lattice} $L$ is a finitely generated $\Z $-submodule of full rank in
Euclidean space $(\R^d,(,))$, so there is some basis 
$(b_1,\ldots , b_d)$ such that $L=\{  \sum _{i=1}^d a_i b_i \mid a_i \in \Z \}$.
The {\em automorphism group } of $L$ is its stabilizer in the 
orthogonal group 
$$\Aut (L) = \{ g\in \GL_d(\R ) \mid g(L) = L  \mbox{ and } (g(x),g(y) ) = (x,y) 
\mbox{ for all } x,y \in \R ^d \} .$$
With respect to a basis of $L$, $\Aut (L)$ is a finite integral matrix group, 
$\Aut (L) \leq \GL_d(\Z )$. 

The {\em dual lattice} is $L^{\#}  := \{ x\in \R^d \mid (x,\ell )\in \Z \mbox{ for all } \ell \in L \}$. We call $L$ {\em integral} if $L\subseteq L^{\#} $ 
and {\em even} if $(\ell,\ell) \in 2\Z $ for all $\ell \in L$. 
The {\em minimum} of $L$ is $\min(L):=\min \{ (\ell,\ell) \mid \ell \in L , \ell \neq 0 \}$.

Even unimodular lattices are
even lattices $L$ with $L=L^{\#}$. They only exist if the 
dimension  is a multiple of 8. 
The theory of modular forms allows to 
bound the minimum of a $d$-dimensional even unimodular lattice by 
$\min (L)  \leq 2 \lfloor \frac{d}{24}  \rfloor +2 $ 
(see for instance \cite{Ebeling})
and lattices achieving
equality are called {\em extremal}.
Very particular lattices are extremal lattice whose dimension is a 
multiple of 24: 
Their sphere packing density realises a local maximum among all 
lattice sphere packings in this dimension and
the vectors of given length in such lattices form 
spherical 11-designs (see e.g. \cite{BachocVenkov}). 
In dimension 24 there is a unique extremal lattice, the 
Leech lattice $\Lambda _{24}$. The Leech lattice is the 
densest lattice in dimension 24, its automorphism group is 
a covering group of the sporadic simple Conway group. 
In dimension 48 one knows 4 extremal lattices (see \cite{aut5}). 
The paper \cite{dim48} shows that there are no other 
48-dimensional extremal lattices that admit a large automorphism.
The aim of the present paper is to establish a similar classification 
in dimension 72. More precisely we prove the following  theorem.

\begin{theorem}
Let $\Gamma $ be an extremal lattice of dimension $72$ that admits an 
automorphism  $\sigma $ such that the minimal polynomial
of $\sigma $ has an irreducible factor of degree $>36$. 
Then $\Gamma = \Gamma _{72}$ is the unique known extremal lattice 
of dimension $72$ constructed in \cite{dim72}.
\end{theorem}

{\bf Acknowledgements.} Many thanks to Steve Donnelly for his 
help in using the class group programs in {\sc Magma} and 
to Markus Kirschmer, who explained me how to use the Mass formulas in
\cite{Hanke}. Part of the programs have been written by Clara Nadenau
in her Master-thesis \cite{Nadenau}.

\section{Some important lattices}

In this paper we denote the $n$-th {\em cyclotomic polynomial} by $\Phi _n$ 
and its degree, the value of the {\em Euler phi function} at $n$,
 by $\varphi (n)$.
To start with let us introduce some important lattices. 
The best studied family of lattices are certainly the {\em root lattices}.
These are even lattices that are generated by vectors $\ell $ of
norm $(\ell,\ell ) = 2$. 
All root lattices are orthogonal sums of the 
indecomposable lattices 
$\A _d$, $\D _d$ ($d\geq 4$), $\E _6$, $\E_ 7$, $\E _8 $ 
(see for instance \cite{Ebeling}).
The next result is certainly well known.

\begin{proposition}
Let $R$ be a root lattice admitting an automorphism $g\in \Aut(R)$ 
of prime order $p>2$ with irreducible minimal polynomial $\Phi _p$.
Then $g$ stabilises each orthogonal summand of $R$. 
If $p=3$ then $R$ is an orthogonal sum of copies of 
$\A_2$, $\D_4$, $\E_6$ or $\E_8$.
For $p=5$ the orthogonal summands are $\A_4 $ and $\E_8$ and 
for $p\geq 7$ the lattice $R$ is an orthogonal sum of copies of $\A _{p-1}$.
\end{proposition} 

\bew
Write $R = R_1^{n_1} \operp  \ldots \operp R_s^{n_s}$ 
with pairwise non isometric indecomposable root lattices $R_i$.
Then $\Aut (R) \cong \Aut(R_1) \wr S_{n_1} \times \ldots \times \Aut(R_s) \wr S_{n_s}$. 
Let $\pi : \Aut(R) \to S_{n_1}\times \ldots \times S_{n_s}$ denote the 
corresponding group homomorphism. 
Then $g\in \ker (\pi )$, because $\Phi_p(g) = 0$, so 
$g$ stabilizes all orthogonal summands of $R$. 
As $\Phi _p$ is irreducible, this shows that $g$ acts on each of the 
summands $R_i$ with minimal polynomial $\Phi _p$. 
Now $\Aut (\A_d) = \langle -1 \rangle \times S_{d+1 }$  contains an
element with minimal polynomial $\Phi _p$ if and only if $S_{d+1}$ has 
a transitive permutation of order $p$, so if and only if $d+1=p$. 
If $d>4$ then $\Aut (\D_{d}) \cong C_2 \wr S_{d} $ contains no
element with minimal polynomial $\Phi _p$ and for 
$\D_4 $,  $\E_6$, $\E_7$ and   $\E_8$ the result follows by 
inspection of the conjugacy classes of elements in the automorphism 
group.
\eb

Closely related to the root lattices are the so called 
{\em Niemeier lattices}, by which we mean the  
 even unimodular lattices of dimension 24.
B. Venkov (\cite[Chapter 18]{SPLAG}) 
has shown that if $L$ is an even unimodular lattice of 
dimension 24, then its {\em root sublattice} $R$,
 which is the sublattice of $L$ generated by the 
vectors of norm 2 in $L$, is either $\{ 0 \}$ (and then 
$L = \Lambda _{24}$ is the Leech lattice, the unique
extremal lattice of dimension 24) or $R$ has full rank. 
Moreover an explicit classification shows that there are 
23 possible $R\neq \{0 \}$ and that $L$ is then 
uniquely determined by $R$.
We will denote the Niemeier lattice $L$ with root sublattice $R$  by $N(R)$.
If $R$ is non empty then $\Aut(N(R)) \leq \Aut (R)$.
As the Leech lattice admits an automorphism of order 7
with irreducible minimal polynomial this shows the 
following corollary.

\begin{corollary} \label{aut7}
Let $L$ be some even unimodular lattice of dimension $24$ admitting 
an automorphism $\sigma $ with minimal polynomial $\Phi _7$.
Then either $L\cong \Lambda _{24}$ or $L\cong N(\A_6^4)$.
\end{corollary}

The fixed lattice of automorphisms of prime order $p$ of even unimodular lattices 
are even $p$-elementary lattices, i.e. even lattices $L$ such that
$pL^{\#} \subseteq L$. 
The possible genera of  such lattices are described in \cite[Theorem 13, Chapter 15]{SPLAG}. The following result holds in greater generality but
we only need it for $p=3$.

\begin{lemma}\label{3elem}
Let $L$ be some even $3$-elementary lattice of dimension $n \leq 24$ such that 
 $\min (L^{\#}) > 2$. Then $L=L^{\#}=\Lambda _{24}$.
\end{lemma}

\bew
As $L$ is an even lattice of odd determinant, the dimension $n=2m$ must
be even. 
By \cite[Theorem 13, Chapter 15]{SPLAG}  there is a 
genus of even 3-elementary lattices of determinant $3^k$ 
and dimension $2m$ if and only if $k\equiv m \pmod{2}$ and this
genus is unique. 
In particular we may always add some orthogonal sum $R$ of 
root lattices $\A _2, \E_6, $ and $\E _8$ so that  $M:=L\operp R$ 
is of dimension 24. By the uniqueness of the 3-elementary genus,
the lattice $M$ is contained in some even unimodular lattice $\Lambda$ 
of dimension 24,
so $$ M = L\operp R \leq \Lambda \leq L^{\#} \operp R^{\#} = M^{\#} .$$
If $\min (L^{\#}) > 2$, then the root system of $\Lambda $ is 
contained in $R^{\#}$, in particular it cannot have full rank.
Venkov's classification of the 24-dimensional even unimodular lattices
shows that this implies that the root system of $\Lambda $ is
empty, so $\Lambda = \Lambda _{24}$ is the Leech lattice, $R=\{0\}$ and
$\Lambda _{24} \subseteq L^{\#} \subseteq \frac{1}{3} \Lambda _{24}$.
By \cite[Theorem 4.1]{Martinet} all non-zero classes of $ \frac{1}{3} \Lambda _{24}/\Lambda _{24} $ contain vectors of norm $\leq \frac{18}{9} = 2$.
So $L^{\#} = \Lambda _{24} = L$.
\eb

Among  the possible fixed lattices of automorphisms of order 5 there 
is one famous lattice that was already studied by J. Tits \cite{Tits} 
and which we call 
${\mathcal T}_{24}$, the {\em Tits lattice}. 
It is a 5-elementary lattice of dimension 24, minimum 8, 
and automorphism group 
$(2.J_2 \circ SL_2(5)) :2$, that is isometric to its 
rescaled dual lattice, ${\mathcal T}_{24}  \cong \sqrt{5} {\mathcal T}_{24}^{\#} $.
It is hence an extremal 5-modular lattice in the sense of 
\cite{Quebbemann}.

\begin{lemma}\label{5elem}
Let $L$ be a $24$-dimensional $5$-elementary even lattice with determinant 
$5^{12}$ and minimum $8$.
\begin{itemize}
\item[(a)] If $L$ has an automorphism $g$ of order $21$ with 
minimal polynomial $\Phi _{21}$, such that  
 the minimal polynomial of the action of 
$g$ on the discriminant group $L^{\#}/L \cong \F_5^{12}$ is also 
$\Phi _{21}$, then $L\cong {\mathcal T}_{24} $.
\item[(b)]
If $L$ has an automorphism $g$ of order $28$ with
minimal polynomial $\Phi_{28}$, such that  
the minimal polynomial of the action of 
$g$ on the discriminant group $L^{\#}/L \cong \F_5^{12}$ is also $\Phi _{28}$,
then $L\cong {\mathcal T}_{24} $.
\end{itemize}
\end{lemma}

\bew
In both cases the minimal polynomial $\mu _g$ is a cyclotomic polynomial and so $L$
can be viewed as some 2-dimensional Hermitian lattice over the ring of
integers in a cyclotomic number field. 
Both parts of the lemma follow by computing the full genus of 
these Hermitian lattices but the techniques are different, due to
the fact that in part (a) the lattices are maximal $\Z[\zeta_{21}]$-lattices 
whereas in part (b), the $\Z[\zeta _{28}]$-lattice $L$ embeds into two
even unimodular $\Z[\zeta_{28}]$-lattices and we may use the
classification of even unimodular 24-dimensional $\Z $-lattices. 
\begin{itemize}
\item[(a)]
Here $\mu_g = f_1 f_2 \in \F_5[x] $ is a product of two irreducible polynomials
modulo 5, where the inverses of the roots of $f_1$ are 
again roots of $f_1$. So as a $\F_5 [g]$-module 
$$L^{\#}/L \cong \F_5[x]/f_1 \oplus \F_5[x]/f_2 \cong S_1\oplus S_2$$ 
is a direct sum of two self-dual simple modules. 
This shows that $L$ is a ${\mathcal A}$-maximal Hermitian lattice, 
where ${\mathcal A}$ is the discriminant of $\Z[\zeta _{21}]$. 
The mass-formula  \cite[Proposition 4.5]{Hanke}, to compute the mass of 
the full genus of Hermitian lattices with $L^{\#}/L \cong S_1\oplus S_2$ to 
be 
$$m_{21} := 2^{-2\cdot 6}\cdot  \frac{48}{63}\cdot  \frac{16}{3}\cdot  2\cdot  (5^3-1)^2 = \frac{1922}{63} .$$
Starting with the lattice $L :={\mathcal T} _{24}$ we use the Kneser neighboring method for both ideals dividing 7 as described in \cite{Schiemann} to 
enumerate the full genus of $L$. 
The class number is $1286$, there are 223 lattice of minimum 4, 
1062 lattice of minimum 6 and a unique lattice of minimum 8, namely 
${\mathcal T} _{24}$.
\item[(b)]
Again 
 $\mu_g = f_1 f_2 \in \F_5[x] $ is a product of two irreducible polynomials
modulo 5, but the inverses of the roots of $f_1$ are 
roots of $f_2$. So as a $\F_5 [g]$-module 
$$L^{\#}/L \cong \F_5[x]/f_1 \oplus \F_5[x]/f_2 \cong S\oplus S^*$$ 
is a direct sum of a simple module $S$ and its dual. So both invariant 
sublattices of $L^{\#}$  that contain $L$ are unimodular. 
In particular for any given lattice $L$ there is a unique 
even unimodular $g$-invariant lattice $M$ containing $L$ such that 
$M/L \cong S$. 
Let $\wp := (f_1(\zeta_{28}) , 5) \unlhd \Z[\zeta _{28}]$ be the 
corresponding prime ideal dividing 5 and view $L$ and $M$ as $\Z[\zeta_{28}]$-lattices
via the action of $g$. 
Then $M/\wp M \cong S\oplus S$ and $L/\wp M$ is one of the $5^6+1=15626$ non-trivial 
submodules of this module. 
\\
There are 2 even unimodular lattices $M$ whose automorphism group
contains an element $a$ of order $7$ with an irreducible minimal polynomial,
$M_1=\Lambda _{24}$, the Leech lattice and $M_2 = N(\A _6^4)$ the Niemeier 
lattice with root system $\A_6^4$. 
The automorphism group of each of these lattices contains a unique element of 
order $28$ with centraliser of order $168$ resp. $196$. 
This again may be verified with the mass formula in \cite{Hanke}, which yields 
$$ 2^{-2\cdot 6}\cdot  \frac{8}{7}\cdot  \frac{416}{21} \cdot 2 = 13/1176 = 1/168 + 1/196 $$ 
as the mass of the genus of $M$. 
So the mass of the genus of $L$ is 
$(5^6+1) \frac{13}{1176} \sim 172.736... $ and one could again use the
Kneser neighboring method to enumerate the genus. 
But here it is much faster to compute the lattices $L$ as sublattices 
of $M_1$ resp. $M_2$. 
There are $2606$ orbits of invariant lattices $L$ of $M_1$ under 
the action of the centraliser of $g$ in $\Aut(M_1)$, 
$703$ of them have minimum 4, 1902 have minimum 6 and a unique sublattice
has minimum 8. 
For $M_2$ one finds 2234 orbits, 2 of minimum 2, 606 of minimum 4, and 1626 of
minimum 6.
\end{itemize}
\eb

To exclude automorphisms of order 25 we classified a certain
genus of 2-dimensional hermitian lattices over $\Z[\zeta_{25}]$:

\begin{lemma}\label{aut25}
Let $Z$ be an even lattice of dimension $40$ that admits an 
automorphism $\sigma \in \Aut(Z)$ with minimal polynomial $ \Phi _{25}$
and such that $(1-\sigma^5) Z^{\#} = Z $. 
Then $\min (Z) \leq 6$.
\end{lemma} 

\bew
To compute the mass of the genus of lattices $Z$ let
$L:=(1-\sigma)^{-2} Z$. Then $L^{\#} = (1-\sigma )^2 Z^{\#} = 
(1-\sigma)^{-1} L $. As there is a bijection between the lattices 
$L$ and $Z$ the mass of the genus of
Hermitian $\Z[\zeta _{25}]$-lattices $L$  and $Z$ coincide. 
Also $L$ is not a maximal lattice, $L^{\#}/L $ is isomorphic to the
hyperbolic plane over $\F_5$ and there are two unimodular $\sigma $-invariant
overlattices $M$ of $L$. 
On the other hand, $M \supseteq L \supseteq (1-\sigma ) M$, so 
$L$ corresponds to one of the 6 one-dimensional subspaces of
$M/ (1-\sigma ) M \cong \F_5^2$.
Using the Venkov trick (as explained in \cite[Proposition 2.4]{BachocNebe})
we hence get for the masses of the genera of Hermitian $\Z[\zeta_{25}]$-lattices
$$ \mbox{mass} (Z) = \mbox{mass} (L) = \frac{6}{2} \mbox{mass} (M) =
\frac{6}{2} \cdot 2^{-2\cdot 10} \cdot \frac{512}{25} \cdot \frac{5825408}{75} \cdot 2 =
\frac{45511}{5000} $$
where we used the mass formula in \cite{Hanke} to compute the mass of the
genus of $M$.
\\
Starting with the lattice $\A_4^{10}$ we enumerate the lattices 
in the genus of $Z$ using Kneser neighboring method until 
we found enough lattices $Z$.
They fall into 50 $\Z $-isometry classes of lattices.
Summing up the reciprocals of the orders of the centralizers 
$C_G(z)$, for representatives of all conjugacy classes of automorphisms 
$z\in G = \Aut_{\Z }(Z) $ such that the minimal polynomial of $z$ is
$\Phi _{25}$ we find almost the full mass up to a summand $3/250$. 
The remaining lattices are constructed ``by hand'': 
There is an obvious even unimodular lattice $M=\E_8^{5}$ that has 
an automorphism of order $25$ with irreducible minimal polynomial. 
Starting with this lattice $M$ we constructed a sublattice $Z$ as illustrated
above. $G=\Aut_{\Z} (Z)$ has 3 conjugacy classes of elements of order
25 with irreducible minimal polynomial and centralisers of order 
$500$, $250$, $250$, respectively. 
The last lattice was found as a sublattice of index 25 of $Z$ invariant
under this group $C_G(z)$ of order $500$ by rescaling the form 
by $(1-z)(1-z^{-1})$. This lattice has the same Hermitian automorphism 
group but already the $\Z $-lattices are not isometric to each other. 
So in total there are 52 $\Z $-isometry classes of lattices $Z$,
which split into 482 isometry classes of 
Hermitian $\Z [\zeta _{25}]$-lattices.
None of the $\Z $-lattices $Z$ has minimum $\geq 8$.
\eb

A remark seems to be in place how to handle high dimensional lattices with
the computer. There is no fast algorithm known that can compute the minimum
of a, say, 72-dimensional lattice from its Gram matrix. 
But there are algorithms to find small vectors in these lattices. 
We never needed to prove extremality of any of the computed lattices,
the only task was to find vectors of norm $\leq 6$.
This was done by a combination of the fast available reduction algorithms
in {\sc Magma}. 

\begin{remark} \label{RED}
To prove that an even  lattice $L$ given by its Gram matrix $F$ 
has minimum $<8$ 
 we performed the following reduction steps in {\sc Magma}:
\\
F:=LLLGram(F:Delta:=0.99,Eta:=0.501); \\
F:=SeysenGram(F); \\
F:=PairReduceGram(F); \\
until the minimum diagonal entry of $F$ becomes $\leq 6$. 
If this minimum is $\geq 8$ after $1000$ loops then we replaced 
$F$ by some equivalent matrix by transforming 
it by some small random element in $\GL_n(\Z )$. 
If $L$ is unimodular then 
 $F\in \GL_n(\Z )$  and we may take 
the equivalent matrix to be $F^{-1} = F^{-1} F F^{-tr}$.
We then repeated $1000$ loops of the reduction steps once more. 
If during these procedures no diagonal entry became smaller than $8$,
there is a fair chance that the minimum of $L$ is $\geq 8$.
\end{remark}

\section{Automorphisms of prime order} 

The notion of the type of an automorphism of a lattice $L$
was introduced in \cite{dim48}. It is motivated by the
analogous notion of a type of an automorphism of a code.

Let $\sigma \in \GL_d(\Q )$ be an element of prime order $p$.
Let $\tilde{F}:=\ker (\sigma -1)$  and $\tilde{Z}:=\im (\sigma -1)$.
Then $\tilde{F}$ is the fixed space of $\sigma $ and the action of
$\sigma  $ on $\tilde{Z}$ gives rise to a vector space structure on $\tilde{Z}$ over the
$p$-th cyclotomic number field $\Q[\zeta _p]$.
In particular $d= f+z(p-1)$, where
$f:=\dim _{\Q} (\tilde{F})$  and $z = \dim _{\Q[\zeta _p]} (\tilde{Z})$.

If $L$ is a $\sigma $-invariant $\Z$-lattice, then
$L$ contains a sublattice  $M$ with
$$ L \geq M= (L\cap \tilde{F}) \oplus (L \cap \tilde{Z}) = : F(\sigma ) \oplus Z (\sigma ) = F\oplus Z \geq pL $$
of finite index $[L:M] = p^s$ where $s\leq \min (f,z )$.

\begin{defn}
The triple $p-(z,f)-s$ is called the {\em type} of the element $\sigma \in \GL(L)$.
\end{defn}

This section is devoted to study the possible types of 
automorphisms of prime order extremal even unimodular lattices of dimension
$24m$.
The arguments are similar to the ones in 
 \cite[Section 4.1]{dim48}, where we obtained contradictions by 
considering the lattices $F$ and $Z$.  Just for 
further reference we recall some results from \cite{dim48} and \cite{aut5}
that will be used in the following arguments without further notice.

\begin{remark} \label{pzfs} 
Let $L$ be an even unimodular lattice and $\sigma \in \Aut (L)$
of prime order $p$ and type $p-(z,f)-s$. 
\begin{itemize}
\item[(a)]
 $s\leq \min (z,f) $ and $s\equiv z \pmod{2}$.
\item[(b)] If $s=0$ then $z$ is even.
\item[(c)] If $s=z$ then $\sqrt{p} Z^{\#} $ is the trace lattice 
of an Hermitian unimodular $\Z[\zeta _p]$ lattice of rank $z$.
\item[(d)] If $s=f$ and $p$ is odd, then  $\frac{1}{\sqrt{p}} F$ is an
even unimodular lattice. In particular $f$ is a multiple of $8$.
\end{itemize}
\end{remark}

The arguments below will be based on the fact that 
both lattices $F$ and $Z$ have minimum $\geq \min (L)$ and 
determinant $\leq p^z$. 
So the Hermite function is
$$\gamma (F) \geq (2m+2) /(p^{z/f}) \mbox{ respectively } 
\gamma (Z) \geq (2m+2)/(p^{1/(p-1)}).$$ 
Already Blichfeldt \cite{Blichfeldt} obtained a good upper bound on
the {\em Hermite constant} $\gamma _d$,
the maximum value of the Hermite function on $d$-dimensional lattices:
$$\gamma _d \leq \frac{2}{\pi} \Gamma (2+\frac{d}{2}) ^{2/d}  = : B(d) .$$
The currently best upper  bounds  for $d\leq 36$ are given
in \cite{Elkies}.

\begin{theorem}\label{allg}
Let $L$ be an extremal even unimodular lattice of dimension $24m$ and
$\sigma \in \Aut(L)$ of prime order $p$ and 
type $p-(z,f)-s$ with $z=1$. Then 
$p=23$, $m=1$, and $L= \Lambda _{24}$ or 
$p=47$, $m=2$, and $L=P_{48q}$. 
\end{theorem} 

\bew
By  Remark \ref{pzfs} 
we get $s\leq \min (z,f) $ and $s \equiv z \pmod{2}$.
So if $z=1$ then $s=1$ and in particular $f\geq 1$.
This excludes that $p=24m+1$ (this also follows from \cite{Bayerirred},
as the characteristic polynomial of $\sigma $ is then the
irreducible polynomial $\Phi_p$).
\\
Assume now that $p=24m-1$, so $f=2$. 
Then $F$ is a 2-dimensional lattice of minimum $2m+2$ and
determinant $p =24m-1$, so 
$$\gamma (F)^2 = \frac{(2m+2)^2}{24m-1} .$$
As the hexagonal lattice is the densest 2-dimensional lattice we have 
$\gamma (F)^2 \leq 4/3$ from which we immediately compute that $m\leq 6$.
If $m=4,5,6$ then $24m-1$ is not a prime. 
The case $m=1$ is clear because the Leech lattice has an automorphism of 
order 23, for $m=2$ the result is \cite[Theorem 5.6]{dim48}.
It remains to consider the case $m=3$, so $p=71$. 
This case has been already treated by Skoruppa \cite{Skoruppa}.
 Then $Z$ is an ideal lattice in $\Q[\zeta _{71}]$ of determinant $71$ and 
minimum 8.  With Magma we checked that all such ideal lattices have minimum
$\leq 6$. 
\\
So from now on we may assume that $p\leq 24m-3$.
\\
First assume that $p\leq 12m-1$. 
Then $Z$ is $(p-1)$-dimensional lattice of determinant $p$ and minimum $\geq 2m+2$.
By Blichfeldt's bound we get 
$$  \frac{2m+2}{p^{1/(p-1)}} \leq \gamma(Z)  \leq B(p-1) = \frac{2}{\pi } 
((\frac{p-1}{2}+1)!)^{2/(p-1)} .$$
Putt $a:=\frac{p-1}{2}$ and take the $a$th power to obtain 
$$(\star ) \ \ (\pi (m+1))^a \leq \sqrt{2a+1} (a+1) a! \leq 
\sqrt{2a+1} (a+1) 
(1+\frac{1}{11a}) \sqrt{2\pi a} \left(\frac{a}{e} \right) ^a $$ 
by Stirling's formula.
Now we use the assumption that $p \leq 12m-1$, so $m\geq \frac{a+1}{6}>\frac{a}{6}$ to conclude that 
$$\left( \frac{e\pi}{6} \right) ^a \leq \sqrt{\pi} (a+1)(2a+1)(1+\frac{1}{11a}) < 6a^2 $$
for $a> 3$. 
As $\frac{e\pi}{6} > \sqrt{2}$ we find 
$2^{a/2} \leq 6 a^2 $ which implies that $a\leq 23$. 
As we know all automorphisms of the Leech lattice and the case $m=2$ is
already treated in \cite{dim48} we may assume that $m\geq 3$. 
But no $1\leq a\leq 23$ satisfies the first inequality in $(\star )$ for
$m=3$. 
\\
Now assume that $12m+1\leq p \leq 24m-3 $.  
Then we consider the lattice $F$ of determinant $p$, 
minimum $\geq 2m+2$, and dimension $2a$ with 
$$2\leq a=(24m-p+1)/2 \leq (24m-(12m+1)+1)/2 = 6m.$$
Blichfeldt's bound gives us 
$$(\star \star) \ \frac{2m+2}{p^{1/(2a)}} \leq \gamma (F) \leq \frac{2}{\pi} ((a+1)!)^{1/a} $$
and as above we take the $a$th power and 
apply Stirling's formula to simplify the right hand side and obtain
$$ \left( \frac{e\pi(m+1)}{a} \right)^a \leq \sqrt{p} \sqrt{2\pi a} (a+1) (1+\frac{1}{11a}) .$$
To bound $\sqrt{p}$ on the right hand side, we divide by 
$\left( \frac{e\pi(m+1)}{a} \right)^{1/2}$ to obtain 
$$ \left( \frac{e\pi(m+1)}{a} \right)^{a-\frac{1}{2}} \leq \sqrt{\frac{p}{m+1}} \sqrt{\frac{2}{e}} a (a+1) (1+\frac{1}{11a}) .$$
As $a\leq 6m \leq 6m+6$ the left hand side is again lower bounded
by $(e\pi/6 )^{a-\frac{1}{2}}$ and hence by $\sqrt{2}^{a-\frac{1}{2}}$. 
For the right hand side (rhs) we use that $p\leq 24 m \leq 25(m+1)$ and
$2/e \leq 1$ to obtain
$\mbox{(rhs)} \leq 5a(a+1)(11a+1)/(11a) \leq 6(a-\frac{1}{2})^2$ whenever $a\geq 12$.
As before we conclude that $a\leq 23$. 
For each value of $a \in \{ 2,\ldots ,23 \}$ 
 the inequality
 $(\star \star)$ can be expressed as the negativity of 
a polynomial in $m$ of degree $2a$ 
(after substituting $p = 24m-2a+1$ and rising to the $2a$-power) for which 
explicit computations show that the unique negative coefficient 
is the one of $m^1$. This gives an easy bound on the maximal 
$m$ for which the polynomial could possibly be positive. 
For $a>1$ this $m$ is always $<4$ and it is immediately
 checked that these polynomials are always positive for $2\leq a\leq 6m+1$. 
\eb

Similar arguments can also be applied to treat other series of parameters, 
but as these are usually immediately ruled out by direct computations 
in any concrete dimension, we omit the heavy calculations. 
Note also that one might not expect such a general theorem for even $z$,
because here $24m=z(p-1)$ and $L=Z$ is always one possibility.

\subsection{Dimension 72}

\begin{theorem}\label{mainp}
Let $\Gamma $ be some extremal even unimodular lattice of dimension $72$ 
and $-I_{72} \neq \sigma \in \Aut(\Gamma )$ be some element of prime order $p$.
Then either $p\in \{37,19,13,7 \}$ and the minimal polynomial of 
$\sigma $ is irreducible (so $\sigma $ has type $p-(72/(p-1),0)-0$) 
or $p\leq 5$ where the possible types are 
\begin{itemize}
\item[p=5]
$5-(z,72-4z)-z$ with $9\leq z\leq 12$,  \\
$5-(13,20)-a$ with $a=13, 11$,  \\
$5-(14,16)-a$ with $a=14,12,10$,  \\
$5-(16,8)-8$, or $5-(18,0)-0$.
\item[p=3]
$3-(z,72-2z)-s$ with $18\leq z\leq 23$ or  \\
 $3-(24,24)-24$ or
 $3-(36,0)-0$.
\item[p=2]
$2-(z,72-z)-s$ with $z=24,32,36,40,48$ 
where $z=s$ if $z=24$ or $z=48$.
\end{itemize}
\end{theorem}

\bew
Let $\sigma \in \Aut(\Gamma )$ be an element of prime order $p$.
By Theorem \ref{allg} we may assume that $z\geq 2$, in particular 
 $p \leq 37$. 
We denote the fixed lattice $F(\sigma )$ by $F$ and $Z(\sigma ) $ by $Z$.
Then $Z$ is an ideal lattice in the $p$-th cyclotomic number field.
\begin{itemize}
\item[$p=37$] By Theorem \ref{allg} we obtain that 
$z=2$.
\item[$p=23,29,31$] Then  $z=2$ and 
in all cases $\gamma (F) > B(72-2(p-1))$, or $z=3$ and $p=23$. 
Then $F$ is a 6-dimensional lattice of determinant $23^3$ and minimum
$8$. So the density of $F$ exceeds the one of 
the root lattice $\E_6$ which is known to be
 the densest lattice of dimension 6. 
\item[$p=19$] Here Blichfeldt's bound leads to the conclusion that 
$z=4$ is the only possibility.
\item[$p=17$] Using Blichfeldt's bound we are left with the possibility
that $z=4$. Then the fixed lattice has dimension 8, determinant $17^4$ and minimum 8.
There is a unique genus of $8$-dimensional 17-elementary even lattices of 
determinant $17^4$, which we enumerated completely using the 
standard procedure implemented in Magma. 
Its class number is 118 and the maximal possible minimum is 6.
\item[$p=13$] Here Blichfeldt's  bound implies that $z=6$.
\item[$p=11$] Using Blichfeldt's bound we are left with the possibility
that $z=6$. Then the fixed lattice is in the genus of the 11-modular 
lattices of dimension 12. By \cite{NebeVenkov} this genus contains no
lattice of minimum 8. 
\item[$p=7$] Blichfeldt's bound allows us to conclude that 
$z\in \{ 9,10,11,12 \}$. 
If $z=11$, then the fixed lattice is a 7-elementary even lattice of dimension 6
and determinant $7^5$. There is a unique genus, this genus has class number 1, 
and the unique lattice in this genus is $\sqrt{7} \A_6^{\#}$ with minimum 6. 
\\
If $z=10$, then the 
fixed lattice is a 12-dimensional, 7-elementary even lattice of 
 determinant $7^{10}$. There is a unique genus, this has class number 12 
and the maximal possible minimum is 6. 
\\
If $z=9$ then the fixed lattice lies in the genus of the 7-modular
lattices of dimension 18. By \cite[Theorem 6.1]{BachocVenkov} this genus 
contains no lattice of minimum 8.
\item[$p=5$] Here we may use Blichfeldt's bound  to obtain 
$z\in \{9,10,11,12,13,14,15,16,18 \}$.
To exclude $z=15$ note that in this case the fixed lattice of $\sigma $ has
dimension 12 and determinant $5^{11}$, more precisely it is in the
genus of $\sqrt{5} (\E_8 \operp \A _4)^{\#} $. 
This genus has class number 2 and both lattices have minimum 4. 
All the other restrictions on the values of $s$ are obtained by 
comparing Blichfeldt's bound with $\gamma (F)$.
Only to exclude the type $5-(12,24)-10$ we need the slightly stronger
bounds in \cite{Elkies}.
\item[$p=3$] The cases $z\leq 14$ are excluded by Blichfeldt's bound.
To exclude the cases $z=15,16,17 $ we look at the lattice $Z$ 
which is (or contains) a unimodular $\Z [\zeta _3]$-lattice. 
So viewed as a $\Z $-lattice this is a $3$-modular $\Z $-lattice 
of dimension $2z =30,32,34$.
The bound in \cite{Quebbemann} shows that the maximal possible minimum of
such a lattice is $6$. 
The cases $25\leq z\leq 35$ are not possible by Lemma \ref{3elem} applied to 
$L= \sqrt{3} F^{\#} $ whence $L^{\#} = \frac{1}{\sqrt{3}} F $ of minimum
$\geq \frac{8}{3} > 2$.
The same lemma also implies that $z=s$ if $z=24$ and then $F =\sqrt{3}\Lambda _{24}$.
\item[$p=2$] 
If $\sigma $ is an automorphism of type $2-(z,f)-s$, then 
 $-\sigma $ has type $2-(f,z)-s$, so we may assume that $f\leq z$. 
The fixed lattice $F$ of $\sigma $ is hence of dimension  $f\leq 36$.
By \cite[Lemma 4.9]{dim48} the lattice $F$ contains $\sqrt{2} U$ for some
unimodular lattice $U$. 
As $\min (F) \geq 8$, we conclude that $\min (U) \geq 4$ and hence 
$f\in \{ 24,32,36 \}$ by the results in 
\cite[Theorem 2]{CoS} and \cite{odd}.
\end{itemize}
\eb

\begin{corollary}
If $p^2 $ divides $|\Aut(\Gamma )|$ then $p\leq 5$. 
\end{corollary}

\bew
Assume that $p>5$. Then  by Theorem \ref{mainp} 
the automorphisms of order $p$ act with
an irreducible minimal polynomial. This shows that $C_p\times C_p$ is not
a subgroup of $\Aut(\Gamma )$. Also as $p$ does not divide $72$ there 
is no element of order $p^2$ in $\Aut(\Gamma )$, which shows the corollary.
\eb

\begin{proposition}
$\Aut(\Gamma )$ does not contain an element of order $25$.
\end{proposition}

\bew
Assume that $\sigma \in \Aut(\Gamma )$ has order $5^2$.
Then $\sigma ^5$ has type $5-(z,f)-s$ so that  $z$ is a multiple of 5. 
By Theorem \ref{mainp} the type of $\sigma ^5$ is $5-(10,32)-10$.
So $Z(\sigma ^5)$ is a lattice as in Lemma \ref{aut25} 
and hence  $\min (Z(\sigma^5)) \leq 6$.
\eb

\begin{remark}
Let $\Gamma _{72}$ be the extremal even unimodular lattice of dimension $72$
constructed in \cite{dim72}. 
Then $\Aut(\Gamma _{72})$ has order $2^83^25^27\cdot 13$
 and contains elements of type
$13-(6,0)-0$, $7-(12,0)-0$, $5-(18,0)-0$, $5-(12,24)-12$ where
the fixed lattice is ${\mathcal T}_{24}$, $3-(24,24)-24$ where
the fixed lattice is $\sqrt{3}\Lambda _{24}$, $2-(48,24)-24$ where
the fixed lattice is $\sqrt{2}\Lambda _{24}$, $2-(36,36)-36$, and 
$2-(24,48)-24$.
\end{remark}

\subsection{Some remarks on dimension 96, 120, 144.} 

The bounds on the density of the fixed lattice $F$ or its orthogonal 
complement $Z$ seem to become more and more restrictive 
in higher dimensions. To illustrate this I give some preliminary 
results in the next three dimensions. 

\begin{theorem}
Let $L$ be an extremal even unimodular lattice of dimension 96, 120, 
or 144 and $\sigma \in \Aut(L)$ some automorphism of prime order 
$p > 7$. 
Then one of
\begin{itemize}
\item $\dim (L) =96$, $p \in \{ 13,17\} $ and $\mu _{\sigma } = \Phi _p$, 
or the type of $\sigma $ is $11-(8,16)-8$ or $13-(7,12)-7$.
\item $\dim (L) =120$, $p\in \{ 11, 13, 31, 61 \} $ and $\mu_{\sigma } = \Phi _p$,
or the type of $\sigma $ is $13-(9,12)-9$.
\item $\dim (L) =144$,  $p\in \{ 13, 19, 37 \} $ and $\mu _{\sigma }  = \Phi _p $.
\end{itemize}
\end{theorem} 

\bew
First assume that $\dim(L) = 96$ and $\min (L) =10$. Combining 
Blichfeldt's bound with the known values for $\gamma _d$ for $d\leq 8$,
excludes all possible types except the ones in the Theorem and 
$11-(9,6)-5$ and $19-(5,6)-5$. 
In these cases a
complete enumeration of the genus of $F(\sigma )$ shows that there is 
no lattice of minimum $\geq 10$.
\\
If $\dim (L) =120$ then Blichfeldt's bound (we 
need the stronger bounds in \cite{Elkies} to exclude the type 
$11-(10,20)-10$) leaves us with the 
possibilities in the Theorem, 
 type $11-(9,10)-9$, and type $17-(7,8)-7$. 
A complete enumeration of the genus of $F(\sigma )$ in the lattes two cases (class number 3) shows that there is no such lattice $F(\sigma )$ of minimum $\geq 12$.
\\
Now assume that $\dim (L) =144$. Then 
Blichfeldt's bound give us with the possibilities in the theorem or
type $29-(5,4)-3$, $13-(11,12)-11$,  $11-(13,14)-13$. 
If $\sigma $ has type $29-(5,4)-3$ then the $F(\sigma )$ is 
one of the two even 29-elementary lattices of dimension 4 and determinant $29^3$.
They have minimum 6 and 4 respectively.
If $\sigma $ has type $13-(11,12)-11$ then the $F(\sigma )$ is 
in the genus of $\sqrt{13}\A_{12}^{\#} $. 
This genus has class number 6 and the minima are $4,4,10,10,12,12$.
In particular no lattice in this genus has minimum $\geq 14 = \min (L)$.
To exclude type $11-(13,14)-13$ we enumerated the genus of the 
even 11-elementary lattices of dimension 14 and determinant $11^{13}$.
The class number is 8 and no lattice in this genus has minimum $\geq 14$. 
\eb

\section{Ideal lattices} 

This section describes the computations to
classify all 72-dimensional extremal even unimodular lattices $L$ that have an 
automorphism $\sigma $ of order $m$ such that 
characteristic polynomial of $\sigma $ is the $m$-th cyclotomic 
polynomial. 
In other words the lattice
 $L$ is an ideal lattice in the $m$-th cyclotomic number field $K = \Q[\zeta _m]$.
By \cite{Bayer} for all these fields there is some positive definite even unimodular ideal lattice. 

\begin{rem} 
Let $\sigma  \in \Aut (L)$ be an automorphism of the lattice $L$ with characteristic polynomial 
$\Phi _m$, the $m$-th cyclotomic polynomial.
Then the action of $\sigma $ on $\Q L$ turns the vector space $\Q L$ into a one-dimensional 
vector space over 
the $m$-th cyclotomic number field $K = \Q[\zeta _m]$ which we identify with $K$.
Then the lattice $L$ is a $\Z [\zeta _m]$-submodule, hence isomorphic to a
fractional  ideal  $J$ in $K$. 
The symmetric positive definite bilinear form $B:L\times L \to \Q $ is $\zeta _m$-invariant,
since $B(\sigma(x) , \sigma(y) ) = B(x,y) $ for all $x,y \in L$. 
In particular  
$$(L,B) \cong (J,b_{\alpha }) \mbox{, where } 
b_{\alpha }(x,y)=\trace _{K/\Q} (\alpha x \overline{y} ) \mbox{, for all } x,y\in J  .$$
Here $\overline{\phantom{x}}$ is the complex conjugation on $K$, the involution 
with fixed field $K^+:= \Q [\zeta _m + \zeta _m^{-1}]$, and $\alpha \in K^+$ 
is totally positive, i.e. $\iota (\alpha ) > 0$ for all embeddings $\iota : K^+ \to \R $. 
\end{rem}

\begin{definition}
Let  $K= \Q[\zeta _m]$ with maximal real subfield $K^{+} := \Q[\zeta _m + \zeta _m^{-1}]$. 
Let $\Cl(K)$ denote the ideal class group of $K$ and $\Cl^+(K^+)$ the ray class group of $K^+$,
where two ideals $I,J $ are equal in $\Cl^+(K^+)$ if there is some totally positive $\alpha \in K^+$
with $I=\alpha J$. Since $x\overline{x} $ is totally positive for any $0\neq x\in K$, 
the norm induces a group homomorphism 
$$ \Cl(K) \to \Cl^+(K^+), \ [J] \mapsto [J\overline{J} \cap K^+]. $$
The {\em positive class group}  $\Cl^+(K)$ is the kernel of this homomorphism and $h^+(K) := | \Cl^+(K) | $. 
\end{definition}

The proof of the next lemma is direct elementary computation
following the lines of \cite{Bayer} and \cite{dim48}.

\begin{lemma}
Let $J$ be a fractional ideal of $K$, $\alpha \in K^+$ totally positive,
and $\Gamma := \Gal(K/\Q) $  the Galois group of $K$ over $\Q $.
\begin{itemize}
\item[(a)] 
Let $$\Delta := \{ x\in K \mid \Tr _{K/\Q } (x \overline{y} ) \in \Z \mbox{ for all } 
y\in \Z [\zeta _m ] \} $$
denote the different of the ring of integers $\Z _K = \Z [\zeta _m]$ of  $K$.
Then the dual lattice of $(J,b_{\alpha }) $ is
$$(J,b_{\alpha })^{\#} = (\overline{J}^{-1} \Delta \alpha ^{-1} ,b_{\alpha }) .$$
So $(J,b_{\alpha })$ is unimodular, if and only if $(J\overline{J})^{-1} \Delta \alpha ^{-1} = \Z_K$.
\item[(b)] $(J,b_{\alpha }) \cong (g(J) , b_{g(\alpha )} ) $ for all 
$g\in \Gamma $.
\item[(c)]  $(J,b_{\alpha } ) \cong (J, b_{u\overline{u} \alpha })$ for 
all $u\in \Z[\zeta_m] ^*$.
\item[(d)]
For $0\neq a\in K$ we have 
$(J,b_{\alpha }) \cong (aJ,b_{\alpha (a\overline{a})^{-1}} ) $. 
\end{itemize} 
\end{lemma}

From this  lemma we  get the following strategy for finding 
representatives of all isometry classes of even unimodular 
ideal lattices:

\begin{rem}\label{reps}
Let $J$ be a fractional ideal of $K$, $\alpha \in K^+$ totally positive,
so that $(J,b_{\alpha })$ is unimodular. 
Then all isometry classes of even unimodular ideal lattices in $K$ are
represented by 
$$(IJ , b_{v \alpha a_I^{-1}}) \mbox{, where } $$
\begin{itemize} 
\item[(a)] $[IJ ] $ runs through a system of representatives of 
the $\Gamma $-orbits on $\Cl^+(K) [J]$,
\item[(b)] $a_I \in K^{+}$ is a totally positive generator of 
$I\overline{I} \cap K^{+} = \Z_{K^+} a_I $, and 
\item[(c)] $v$ runs through a set of 
representatives of $U^+/ N$, where
$$U^+:=(\Z_{K^+}^*)_{>> 0}  = 
\{ v\in  \Z_{K^+}^* \mid v \mbox{ is totally positive } \} $$
is the group of totally positive units in $\Z_{K^+}$ and 
 $$N:=N(\Z _{K}^*) := \{ u\overline u \mid u\in \Z_{K}^* \} $$
the subgroup of norms of units in $\Z_{K}$.
\end{itemize} 
\end{rem}

\knubbel 
To compute a system of representatives as in Remark \ref{reps} (c) 
we use \cite[Theorem 8.3]{Washington}, that gives a system 
$\underline{u}:=(u_1,\ldots , u_{r-1})$ of multiplicatively
independent units in $\Z_{K^+}^*$ and the index $i = [
\Z_{K^+}^*:\langle \underline{u} \rangle]$. 
To find the totally positive units $U^+$ 
(modulo the squares)  we put $u_r :=-1$ and compute a matrix 
$R(\underline{u}) \in \F_2^{r\times r} $ with 
$$R(\underline{u})_{i,j} = \left\{ \begin{array}{ll} 1 & \mbox{ if } \epsilon_j(u_i)  < 0 \\
 0 & \mbox{ if } \epsilon_j(u_i)  > 0  \end{array} \right. $$
for $i=1,\ldots , r $ and all real embeddings 
$\epsilon _1,\ldots , \epsilon _r $
of $K^+$. 
Any $s$ in the kernel of $R(\underline{u})$ corresponds to a totally positive
unit $\prod _{i=1}^r u_i^{s_i}$ which is a square in $\langle u_1,\ldots , u_r \rangle $, if and only if $s=0$. 
If $i=[\Z_{K^+}^*: \langle u_1,\ldots , u_r \rangle ]$ is odd, then 
the kernel of $R(\underline{u})$ gives a system of representatives of the totally
positive units modulo $(\Z _{K^+}^*)^2$. 
If not, then one of these totally positive elements 
$u = \prod _{i=1}^r u_i^{s_i}$  is a square in $K^+$ and
we find its square root $v\in \Z _{K^+}^*$  by factoring $X^2-u$ 
in $K^+[X]$ and replace one of the $u_i$ with $s_i=1$
by $v$ until we find a subgroup 
$V :=\langle u_1,\ldots , u_r \rangle $ with $[\Z_{K^+}^* : V]$ odd. 
Then $U^+ /(\Z _{K^+}^*)^2 \cong V_{>>0} / V^2$. 

\knubbel
By \cite[Satz 14 and Satz 27]{Hasse} $[N : (\Z _{K^+}^*)^2 ] = 2$
whenever $m$ is not a prime power. 
To find $N$ without computing $\Z_K^*$ we use the 
Galois action on $U^+ /(\Z _{K^+}^*)^2 \cong V_{>>0} / V^2$.
This group is an $\F_2 \Gamma $ module containing the trivial
$\F_2\Gamma $ submodule
$N/(\Z _{K^+}^*)^2$. 
In all cases this $\F_2\Gamma $ module turned out to have a  unique trivial submodule
which is hence the image of $N$.

\knubbel
To obtain suitable representatives of the Galois orbits on the 
ideal classes as in Remark \ref{reps} (a) we usually 
choose a totally decomposed small prime $p$ and 
compute the permutation action of $\Gamma $ on the 
set of prime ideals $\{ P_1,\ldots , P_{72} \}$ of $\Z_K$ 
that contain $p$. In most cases it was easy to 
find a small prime $p$ such that 
$$\langle [P_1],\ldots , [P_{72}] \rangle = \Cl^+ (K) $$
and to compute the explicit $\Z \Gamma $-module epimorphism 
$$f:\Z^{72} \to \Cl^+(K) , v:=(v_1,\ldots, v_{72}) \mapsto \prod _{i=1}^{72}  [P_i ]^{v_i} $$
If not we take several primes $p$ and consider the direct product of
permutation modules. 
We successively compute small elements $v\in \Z^{72}$ 
such that $\bigcup _{v} \{ f(w ) \mid w\in v^{\Gamma }  \} = \Cl^+(K)$. 
Of course, if $J\neq 1$, then we need to compute the action of $\Gamma $ 
on the elements $[J \prod _{i=1}^{72} P_i^{v_i}]$. 
In our cases we could choose $[J]$ to be a class of order $4$, so
$[J]^{\gamma } = [J]^{\pm 1}$ for all $\gamma $ in $\Gamma $.
As 
$[J \prod _{i=1}^{72} P_i^{v_i}] ^{\gamma } = [J]^{\gamma } 
[\prod _{i=1}^{72} P_i^{v_i} ]^{\gamma }$
this can be controlled by adding an additional block diagonal entry 
$\pm 1$ to the matrices representing the $\Gamma $-action on 
$\Z^{72}$ and replacing $f$ accordingly.

\knubbel
The following table summarizes the computations that 
led to the theorem below. 
$K= \Q[\zeta _m]$ with $m\not\equiv 2 \pmod{4} $ and $\varphi (m) = 72$. 
Then $h_K$ displays the structure of the class group 
 as computed with {\sc Magma}. The class number 
coincides with the one given in \cite{Washington}. In all cases 
$\Cl^+(K) = \Cl(K)^2$. The third row gives the index of the 
norms in the group of totally positive units and line four 
an ideal class $[J]$ 
that supports even unimodular positive definite lattices. 
We then list the primes $p$ for which the prime ideals of
$\Z_K$ that divide $p$ generate $\Cl^+(K)$ followed by the
number of Galois orbits on $\Cl^+(K)$.
$$
\begin{array}{|l|c|c|c|c|c|} 
\hline
m & 91 & 95 & 111 & 117 & 135  \\
h_K & 4\times 13468 & 107692  & 19\times 25308 & 3\times 3\times 9 \times 1638 & 75961 \\
|U^+/N | &  4 & 2 & 2 & 2 & 1 \\
 J & 1 & Q_{11}  
 & Q_{1999} 
&
1 
&
1 \\
p & 547^2 & 191 &  223 & 1171 & 271 \\
\# orbs & 202 & 751 & 3433 & 942 & 1086 \\
\hline
m & 148 & 152 & 216  & 228  & 252 \\
h_K  & 4827501 & 19\times 171\times 513 & 1714617  & 238203 & 7\times 28\times 364 \\
|U^+/N | &  1 & 1 & 1 & 1 & 4 \\
  J & 1 & 1 & 1 &1 & 1 \\
p &  593 & 457 & 433 & 229 &  2,5,71 \\
\# orbs & 68891^* & 24568^* & 25337^* & 3878 & 338  \\
\hline
\end{array}
$$
Here $Q_{11}$ resp. $Q_{1999}$ denotes a suitable product of three 
prime ideals in $\Z_K$ that divide $11$ resp. $1999$ so that 
the classes $[Q_{11}]$ and $[Q_{1999}]$ have order 4.
\\
For $m=91$ the entry $p=547^2$ means that we took the squares of 
the prime ideals dividing $547$, as these generate $\Cl^+(K) = \Cl(K)^2$.
\\
Some of the numbers of orbits are marked with a $*$: \\
For $m=148$ our random procedure generating small elements  
$v\in \Z ^{72}$ such that $\bigcup _v f(v^\Gamma ) = \Cl ^+(K)$ 
did not find certain elements of order 9 and 3. 
So additionally to the 68891 orbit representatives we considered all
elements in  $C_{9} \cong \langle [P_{37}] \rangle \leq \Cl ^+(K) $.
\\
Similarly for $m=152$ the $24568$ Galois orbits did not 
represent some elements of order 9 in the class group.
We additionally considered all elements in the Sylow 3-subgroup 
$(\cong C_9\times C_{27})$ of the class group.
\\
Finally for $m=216$ we also failed to find the six elements of order 9
in the class group using the random search for small representatives. 
These are found in the subgroup 
 $\langle [P_{3}] \rangle \cong C_9 $
generated by 
the prime ideals dividing 3.

Using the strategy explained before we find the following result.

\begin{theorem}\label{ideal} 
Let $L$ be an extremal even unimodular lattice such that $\Aut (L)$ contains some 
element $\sigma $ with irreducible characteristic polynomial.
Then the order of $\sigma $ is $91$ or $182$  and $L=(I,b_{\alpha }) \cong \Gamma _{72}$ for all six ideal
classes $[I]\in \Cl(\Q[\zeta_{91}])$ of order $7$ and 
a unique $\alpha $ modulo $N(\Z [\zeta _{91}]^*)$.
\end{theorem}

To prove this theorem we computed representatives of all isomorphism classes 
of ideal lattices in the $m$-th cyclotomic number field as explained above 
using the class group algorithm available in {\sc Magma}. 
We then performed the reduction procedure described in Remark \ref{RED}.
The lattices that survived with minimum $\geq 8$ correspond to 
ideal lattices $(I,b_{\alpha })$ for 
the six ideal classes $[I]$ of order 7 in $\Cl (\Q [\zeta _{91}])$ with 
a  unique $\alpha $ modulo squares, so these are all in the same
Galois orbit and hence isometric.
As the lattice $\Gamma _{72}$ from \cite{dim72} has an automorphism
of order 91 with irreducible characteristic polynomial, this 
ideal lattice has to be the extremal even unimodular lattice $\Gamma _{72}$.

\section{Large automorphisms.}

Similar as ideal lattices we can classify lattices that 
admit a large automorphism.

\begin{definition} 
Let $\sigma \in \GL_d(\Z )$. Then $\sigma $ is called {\em large},
if the minimal polynomial of $\sigma $ has an irreducible
factor of degree $>d/2$. 
\end{definition} 

If  $\sigma \in \Aut(L)$ is a large automorphism 
then the irreducible factor of the minimal polynomial is some cyclotomic
polynomial $\Phi _n$ with $\varphi (n) > d/2$ and 
the lattice $L$ contains a sublattice 
$Z \operp F$ of finite index
where  $Z$ is an ideal lattice in the $n$-th 
cyclotomic number field. 
The lattice $L$ is an ideal lattice, if $F=\{ 0 \}$ or equivalently 
$\varphi(n) = d$.

\begin{remark}
If $\sigma \in \Aut(L) $ is a large automorphism then the order $n$ of
$\sigma $ satisfies $\varphi (n) > d/2$. 
\end{remark}

This section describes the computations that lead to the following result:

\begin{theorem}\label{subideal}
Let $\Gamma $ be an extremal even unimodular lattice of dimension $72$
and assume that there is some large automorphism
$\sigma \in \Aut(\Gamma )$. Then the order of $\sigma $ is 
either $91$, $182$ or $168$ and 
$\Gamma = \Gamma _{72} $. 
\end{theorem}

Note that $\Aut (\Gamma _{72})$ has four conjugacy classes of elements of 
order 168. All these elements $g$ have minimal polynomial 
$\Phi _{168} \Phi_{56} $ of degree 72. 
So if $F$ is the fixed lattice of $g^{56}$ then $F\cong \sqrt{3}\Lambda _{24}$
is a rescaled Leech lattice. 
For each of the two conjugacy classes $[h]$ of elements
of order 56 in $\Aut(\Lambda _{24})$ exactly 
two conjugacy classes of elements $g\in \Aut (\Gamma _{72})$ 
satisfy $g_{|F } \in [h]$. 

During the whole section we keep the assumptions of Theorem \ref{subideal}.
As we already dealt with ideal lattices in the previous section
we assume that the characteristic polynomial of $\sigma $ is not irreducible.
In particular $\Gamma $ contains a sublattice $Z\operp F$,
so that $Z$ is an ideal lattice in the $n$-th cyclotomic number field
and $F$ is the fixed lattice of $\sigma ^{n/p}$ for some prime divisor 
$p$ of $n$. 

From Theorem \ref{mainp} we hence obtain the following list of 
possible element orders $n$ such that $n\not\equiv 2 \pmod{4}$ and 
$36< \varphi (n) < 72 $:
\\
$\varphi (n) = 64 $. Then $n=2^7$, $n=2^5 5$, $n=2^6 3 $, or $n=2^4 3 \cdot 5$. 
\\
$\varphi (n) = 54$ and $n=3^4$.  \\
$\varphi(n) = 48$. Then the following table displays the possible values 
for $n$ together with the structure $h_K$ of the class group of $K=\Q[\zeta _n]$ 
and the positive class group $h_k^+$: \\
\begin{tabular}{|l|c|c|c|c|c|c|c|c|c|}
\hline
n & 65 & 104 & 105 & 112 & 140
 & 144 & 156 & 168  & 180   \\
$h_K$ & $2\times 2 \times 4 \times 4$ &$ 3\times 117$  & $13$ & $3\times 156$ & $39$
  & $13 \times 39$ & $156$ & $84$  & $5\times 15$  \\
$h^+_K$ & $2\times 2$ & $3\times 117$ &  $13$ & $6\times 39$ & $39$
 & $13 \times 39$ & $78$ & $42$  & $5\times 15$  \\
$|U_+/N|$ &  16 & 1 & 1 & 2 & 1 & 1 & 2 & 2 &  1 \\
\hline
\end{tabular}
\\
$U_+/N$ is the group of totally positive units in $\Z[\zeta_n + \zeta_n^{-1}]^*$
modulo the group of norms of units in $\Z[\zeta_n]$

\begin{lemma}
$n\neq 3^4$.
\end{lemma} 

\bew
Otherwise let $\sigma \in \Aut(\Gamma )$ be some element of order $3^4$.
Then $\alpha := \sigma ^{27}$ has type $3-(27,18)-s$ which contradicts 
Theorem \ref{mainp}.
\eb

\begin{lemma} 
$\varphi (n) \neq 64 $.
\end{lemma} 

\bew
In this case $F$ has dimension 8. By Theorem \ref{mainp} the only possibility
is that $F$ is the fixed lattice of some element of order $5$, and
then $F= \sqrt{5} E_8$. This lattice has no automorphism of order
$2^5$ or $2^43$. 
\eb

\begin{lemma}
$n\neq 65 = 5\cdot 13$.
\end{lemma} 

\bew
Assume that $n=65$, then $F$ is the fixed lattice of an element of
order 5, so $\dim (F) =24$, $\det (F) = 5^{12}$ and $Z$ is an ideal 
lattice of determinant $5^{12}$.
There is no such positive definite ideal lattice $Z$. 
\eb

\begin{lemma}
$n\neq 104 = 2^3 13$.
\end{lemma}

\bew
If $n=104$, then $F$ is the fixed lattice of an element of order 2 and hence 
$F \cong \sqrt{2} \Lambda _{24}$ and $Z$ is a positive definite ideal lattice
of determinant $2^{24}$.
All ideal classes support a unique such lattice $Z$, but only for those classes 
in a certain cyclic subgroup of order 9 this lattice $Z$ has 
minimum $\geq 8$.  These 9 classes fall into 3 Galois orbits.
In all cases $\zeta _{104}$ acts like a primitive $52$nd root of unity on
$Z^{\#}/Z \cong \F_2^{24}$. 
The automorphism group $\Aut(\Lambda _{24})$ of the Leech lattice has 
a unique conjugacy class of elements $g$ of order 52, so the extremal even unimodular 
72-dimensional lattices $\Gamma $ with an automorphism $\sigma $ of order 104
are maximal isotropic $(\zeta_{104},g)$-invariant submodules of 
$$ X:= (Z^{\#} \operp \frac{1}{\sqrt{2}} \Lambda _{24}) / (Z \operp \sqrt{2}\Lambda _{24})  \cong \F_{2}^{48} $$
The module $X$ has $2^{12}-1$ invariant minimal submodules $M$ that fall 
into $317$ orbits under the action of the centralizer $\langle (1,g) \rangle $. 
For representatives $M$ of these orbits we computed the maximal isotropic
invariant submodules of $M^{\perp }/M$ 
and the corresponding 72-dimensional unimodular lattices.
In all these lattices our reduction algorithm \ref{RED} 
finds vectors of norm $<8$, so none of these lattices is extremal.
\eb

\begin{lemma}
$n\neq 105 = 3 \cdot 5 \cdot 7$.
\end{lemma}

\bew
If $n=105$, then $F$ is the fixed lattice of an element of order $p$ 
for $p\in \{ 3,5 \}$. Moreover $\sigma $ acts on $F^{\#}/F $ as an 
element of order $\frac{105}{p} \in \{ 35, 21 \}$. 
Hence $F\cong \sqrt{3} \Lambda _{24}$ if $p=3$ and by Lemma \ref{5elem} (a) 
$F\cong {\mathcal T} _{24}$ if $p=5$.
So 
$Z$ is a positive definite ideal lattice
of determinant $3^{24}$ or $5^{12}$ and minimum $\geq 8$.
There are two Galois orbits on the ideal classes of $\Z[\zeta _{105}]$, 
both support a unique positive definite lattice of determinant 
$3^{24}$ and $5^{12}$. 
For $p=3$ both lattices have minimum 8, for $p=5$ only the principal 
ideals support lattices of minimum 8. 
\\
\underline{p=3}: 
The automorphism group $\Aut(\Lambda _{24})$ of the Leech lattice has 
a unique conjugacy class of elements $g$ of order 35, so the extremal even unimodular 
72-dimensional lattices $\Gamma $ with an automorphism $\sigma $ of order 105
correspond to maximal isotropic $(\zeta_{105},g)$-invariant submodules $M$ of 
$$ X:= (Z^{\#} \operp \frac{1}{\sqrt{3}} \Lambda _{24}) / (Z \operp \sqrt{3}\Lambda _{24})  \cong \F_{3}^{48} \cong S^2 \oplus (S^*)^2$$
as a $(\zeta_{105},g)$-module, where $S^*$ is the dual module of the simple 
module $S\cong \F_{3}^{12}$. 
So $M = Y \oplus Y^{\perp }$ where $Y\leq S^2$ and $Y^{\perp }$ is its 
annihilator in $(S^*)^2$. 
There are $3^{12}+3=531444$ submodules $Y\leq S^2$ which fall into 
$15187$ orbits 
 under the action of the centralizer $\langle (1,g) \rangle $. 
For representatives $Y$ of these orbits we compute the unimodular lattice $\Gamma $ 
as the full preimage of $Y\oplus Y^{\perp }$.
For both ideals no such lattice has minimum $\geq 8$.
\\
\underline{p=5}:
The automorphism group $\Aut({\mathcal T} _{24})$ of the Tits lattice has 
a unique conjugacy class of elements $g$ of order 21, so the extremal even unimodular 
72-dimensional lattices $\Gamma $ with an automorphism $\sigma $ of order 105
correspond to maximal isotropic $(\zeta_{105},g)$-invariant submodules $M$ of
$$ X:= (Z^{\#} \operp {\mathcal T} _{24}^{\#} ) / (Z \operp {\mathcal T} _{24})  \cong \F_{5}^{24} .$$
The  module $X$ has $2\cdot (5^6+1)=31252$ minimal submodules 
falling into 1492
 orbits 
 under the action of the centralizer $\langle (1,g) \rangle $.
For all representatives we compute the maximal isotropic invariant 
overmodules and the unimodular lattices as full preimages.
None of them has minimum $\geq 8$.
\eb

\begin{lemma}
$n\neq 112=2^47$.
\end{lemma}

\bew
If $n=112$, then $F$ is the fixed lattice of an element of order 2,
so $F\cong \sqrt{2}\Lambda _{24}$ and $Z$ is an ideal lattice of 
determinant $2^{24}$. 
The Galois group acts on the ideal class group with 16 orbits and
all ideals support a unique positive definite lattice of determinant $2^{24}$.
None of these lattices has minimum $\geq 8$.
\eb

\begin{lemma}
$n\neq 140 = 2^2 \cdot 5 \cdot 7$.
\end{lemma}

\bew
If $n=140$, then $F$ is the fixed lattice of an element of order $p$ 
for $p\in \{ 2,5 \}$. Moreover $\sigma $ acts on $F^{\#}/F $ as an 
element of order $35$ resp. $28$. 
Hence $F\cong \sqrt{2} \Lambda _{24}$ if $p=2$ and by Lemma \ref{5elem} (b) 
$F\cong {\mathcal T} _{24}$ if $p=5$.
So 
$Z$ is a positive definite ideal lattice
of determinant $2^{24}$ or $5^{12}$ and minimum $\geq 8$.
There are 4 Galois orbits on the ideal classes of $\Z[\zeta _{140}]$, 
all four support a unique positive definite lattice of determinant 
$2^{24}$ and $5^{12}$. 
For $p=2$ the lattices that have minimum 8 correspond to the 
ideal classes of order 3 and 1, for $p=5$ only the principal 
ideals support lattices of minimum 8. 
\\
\underline{p=2}: 
The automorphism group $\Aut(\Lambda _{24})$ of the Leech lattice has 
a unique conjugacy class of elements $g$ of order 35, so the extremal even unimodular 
72-dimensional lattices $\Gamma $ with an automorphism $\sigma $ of order 140
correspond to maximal isotropic $(\zeta_{140},g)$-invariant submodules $M$ of 
$$ X:= (Z^{\#} \operp \frac{1}{\sqrt{2}} \Lambda _{24}) / (Z \operp \sqrt{2}\Lambda _{24})  \cong \F_{2}^{48} \cong S^2 \oplus (S^*)^2$$
as a $(\zeta_{140},g)$-module, where $S^*$ is the dual module of the simple 
module $S\cong \F_{2}^{12}$. 
So $M = Y \oplus Y^{\perp }$ where $Y\leq S^2$ and $Y^{\perp }$ is its 
annihilator in $(S^*)^2$. 
There are $2^{12}+3=4099$ submodules $Y\leq S^2$ which fall into 
$121$ orbits 
 under the action of the centralizer $\langle (1,g) \rangle $. 
For representatives $Y$ of these orbits we compute the unimodular lattice $\Gamma $ 
as the full preimage of $Y\oplus Y^{\perp }$.
For both ideals no such lattice has minimum $\geq 8$.
\\
\underline{p=5}:
The automorphism group $\Aut({\mathcal T} _{24})$ of the Tits lattice has 
a unique conjugacy class of elements $g$ of order 28, so the extremal even unimodular 
72-dimensional lattices $\Gamma $ with an automorphism $\sigma $ of order 140
correspond to maximal isotropic $(\zeta_{140},g)$-invariant submodules $M$ of
$$ X:= (Z^{\#} \operp {\mathcal T} _{24}^{\#} ) / (Z \operp {\mathcal T} _{24})  \cong \F_{5}^{24}  \cong S^2 \oplus (S^*)^2$$
as a $(\zeta_{140},g)$-module, where $S^*$ is the dual module of the simple 
module $S\cong \F_{2}^{12}$. 
So $M = Y \oplus Y^{\perp }$ where $Y\leq S^2$ and $Y^{\perp }$ is its 
annihilator in $(S^*)^2$.
There are $5^{6}+3=15628$ submodules $Y\leq S^2$ which fall into  
$ 562$ orbits 
 under the action of the centralizer $\langle (1,g) \rangle $. 
For representatives $Y$ of these orbits we compute the unimodular lattice $\Gamma $
as the full preimage of $Y\oplus Y^{\perp }$.
No such lattice has minimum $\geq 8$.
\eb

\begin{lemma}
$n\neq 144 = 2^4 3$.
\end{lemma}

\bew
If $n=144$, then $F$ is the fixed lattice of an element of order $2$  or
$3$. In particular $Z$ is an ideal lattice of determinant $2^{24}$ 
or $3^{24}$ in $\Q[\zeta _{144}]$.
Explicit computations show that $\zeta _{144}$ acts as an element of order
$36$ (with minimal polynomial $(x^6+x^3+1)^4 \in \F_2[x]$) (for $p=2$) 
respectively an element of order 
$48$ (with minimal polynomial $(x^4+x^2+2)^3(x^4+2x^2+2)^3 \in \F_3[x] $) (for $p=3$) 
on $Z^{\#}/Z$. 
This leads to a contradiction, as $\Aut(\Lambda _{24})$ has no such elements acting
this way on $\Lambda _{24}/p\Lambda _{24}$.
\eb

\begin{lemma}
$n\neq 156 = 2^2 3\cdot 13$.
\end{lemma}

\bew
Assume that $n=156$. Then $F$ is the fixed lattice of an element of order 2 or 3,
so $Z$ is an ideal lattice of determinant $2^{24}$ resp. $3^{24}$. 
The Galois group of $\Q[\zeta_{156}]$ has eight orbits on $\Cl ^+(\Q[\zeta_{156}] ) \cong 
C_{78}$ according to the eight different element orders in this group. 
On all ideals $I$ there is some positive definite even unimodular form $B_{\alpha _I }$. 
\\
\underline{p=3:} 
In $\Z[\zeta _{156} + \zeta _{156}^{-1} ]$ the ideal $(3) = (p_3)^2 (p_3')^2 $
with $p_3,p_3' \in \Z[\zeta _{156} + \zeta _{156}^{-1} ]$ so that 
$p_3p_3'$ is totally positive. 
So the possible ideal lattices $Z$ are 
$Z_1:=(I,B_{p_3^2 \alpha  _I} )$, 
$Z_2:=(I,B_{(p_3')^2 \alpha _I } )$, and
$Z_3:=(I,B_{p_3p_3' \alpha _I } )$. 
On $Z_1^{\#}/Z_1$ and also on $Z_2^{\#} / Z_2$ the element $\zeta_{156}$ 
acts as an element of order $156$. Such an element 
does not exist in $\Aut (\Lambda _{24})$. 
So we need to consider the ideal lattices 
$(I,B_{p_3p_3'\alpha_I}) $ and $(I,B_{up_3p_3'\alpha_I})$ for 
all 8 representatives $[I]$ of the Galois orbits and some unit $u$
that is not a norm. 
We find 5 such ideal lattices that have minimum $\geq 8$: 
For the elements $[I]$ of order $39$ in the class group and one unit, 
the elements of order $26$ and all units and the elements of order 1,
again both units. 
The Leech lattice has a unique conjugacy classes of automorphisms of order 52.
As already described before we compute the invariant unimodular overlattices
of $Z\perp \sqrt{3} \Lambda _{24}$ for all five lattices $Z$.
None of them has minimum 8.
\\
\underline{p=2:} 
Here we find on each ideal class $2 = |U^+/N|$ lattices $Z$ of determinant $2^{24}$,
only for the principal ideals these lattices have minimum $\geq 8$. 
The Leech lattice has two conjugacy classes of automorphisms of order 39, 
represented by, say, $g$ and $h$. 
As already described before we compute the unimodular overlattices 
of $Z\perp \sqrt{2} \Lambda _{24}$ for both lattices $Z$ invariant 
under $(\zeta_{156},g)$ respectively $(\zeta_{156},h)$.
None of them has minimum 8. 
\eb

\begin{lemma}
If 
$n=168 = 2^3 3 \cdot 7$ then $\Gamma = \Gamma _{72}$.
\end{lemma}

\bew
Assume that $n=168$. Then $F$ is the fixed lattice of an element of order 2 or 3,
so $Z$ is an ideal lattice of determinant $2^{24}$ resp. $3^{24}$.
The Galois group of $\Q[\zeta_{168}]$ has eight orbits on $\Cl ^+(\Q[\zeta_{168}] ) \cong 
C_{42}$ according to the eight different element orders in this group.
On all ideals $I$ there is some positive definite even unimodular form $B_{\alpha _I }$.
\\
\underline{p=3:}
In $\Z[\zeta _{168} + \zeta _{168}^{-1} ]$ the ideal $(3) = (p_3)^2 (p_3')^2 $
with $p_3,p_3' \in \Z[\zeta _{168} + \zeta _{168}^{-1} ]$ so that
$p_3p_3'$ is totally positive.
So the possible ideal lattices $Z$ are
$Z_1:=(I,B_{p_3^2 \alpha  _I} )$,
$Z_2:=(I,B_{(p_3')^2 \alpha _I } )$, and
$Z_3:=(I,B_{p_3p_3' \alpha _I } )$.
On $Z_1^{\#}/Z_1$ and also on $Z_2^{\#} / Z_2$ the element $\zeta_{168}$
acts as an element of order $168$. Such an element
does not exist in $\Aut (\Lambda _{24})$.
So we need to consider the ideal lattices
$(I,B_{p_3p_3'\alpha_I}) $ and $(I,B_{up_3p_3'\alpha_I})$ for
all 8 representatives $[I]$ of the Galois orbits and some unit $u \in U^+\setminus N$.
We find 9 such ideal lattices that have minimum $\geq 8$:
For the elements $[I]$ of order $2$, $6$, and $14$ in the class group and both units,
the elements of order $21$ and the unit 1, and for the 
elements of order 7 and 1 and the unit $u$.
The Leech lattice has  two conjugacy classes of automorphisms of order 56.
As already described before we compute the invariant unimodular overlattices
of $Z\perp \sqrt{3} \Lambda _{24}$ for all nine lattices $Z$ and representatives 
of both conjugacy classes. Only for the lattices $Z$ where the underlying ideal is principal
the algorithm \ref{RED} failed to find vectors of norm $<8$ in exactly 2 
even unimodular overlattices of $Z\perp \sqrt{3} \Lambda _{24}$ (for both conjugacy 
classes of elements of order 56 in $\Aut(\Lambda _{24})$). 
As $\Aut (\Gamma _{72})$ has four conjugacy classes of elements of order 168 
and these have minimal polynomial 
$\Phi _{168} \Phi_{56} $ these four lattices are indeed isometric to 
$\Gamma _{72}$.
\\
\underline{p=2:}
Here we find on each ideal class $2 = |U^+/N|$ lattices $Z$ of determinant $2^{24}$,
only for the classes of order $1$ and $7$ and the unit 1, 
and the classes of order $2$ and $6$ and both units  these lattices have minimum $\geq 8$.
The Leech lattice has a unique conjugacy class $[g]$ of automorphisms of order 84,
which act as elements of order $42$ on $\Lambda_{24}/2\Lambda _{24}$.
For all six lattices $Z$ invariant
under $(\zeta_{168},g)$ we compute the unimodular overlattices
of $Z\perp \sqrt{2} \Lambda _{24}$.
None of them has minimum 8.
\eb

\begin{lemma}
$n\neq 180 = 2^2 3^2 5$.
\end{lemma}

\bew
If $n=180$, then $F$ is the fixed lattice of an element of order $2$, $3$,  or
$5$. In particular $Z$ is an ideal lattice of determinant $2^{24}$, $3^{24}$, 
or $5^{12}$ in
$\Q[\zeta _{180}]$.
Explicit computations show that $\zeta _{180}$ acts as an element of order
$45$ (for $p=2$)
respectively an element of order
$60$ with minimal polynomial $(x^4+x^3+2x+1)^3(x^4+2x^3+x+1)^3 \in \F_3[x] $ (for $p=3$)
on $Z^{\#}/Z$.
This leads to a contradiction, as $\Aut(\Lambda _{24})$ has no such elements acting
this way on $\Lambda _{24}/p\Lambda _{24}$.
For $p=5$ we compute the minima of the suitable ideal lattices of determinant $5^{12}$
to see that there is no such ideal lattice with minimum $\geq 8$.
\eb

\end{document}